\documentclass[a4paper]{article}

\usepackage{amsmath,amsthm}
\usepackage{cite}
\usepackage{graphicx}
\usepackage{url}
\usepackage{hyperref}

\newtheorem{theorem}{Theorem}[section]
\newtheorem{lemma}[theorem]{Lemma}
\newtheorem{proposition}[theorem]{Proposition}
\newtheorem{corollary}[theorem]{Corollary}
\newtheorem{definition}[theorem]{Definition}
\newtheorem{remark}[theorem]{Remark}

\DeclareMathOperator{\sgn}{sgn}
\DeclareMathOperator{\Sgn}{Sgn}
\DeclareMathOperator{\Span}{span}
\DeclareMathOperator{\range}{range}
\DeclareMathOperator{\dom}{dom}

\DeclareMathOperator*{\argmin}{argmin}

\newcommand{\scp}[2]{\left.\left\langle #1\vphantom{#2}\right|#2\right\rangle}
\newcommand{\norm}[1]{\left\|#1\right\|}

\newcommand{\R}{\mathbf{R}}
\newcommand{\N}{\mathbf{N}}
\newcommand{\bigO}{\mathcal{O}}
\newcommand{\alphaeff}{\alpha_\textup{eff}}

\newcommand{\lpspace}[1]{\ell^{#1}}
\newcommand{\logor}{\vee}
\newcommand{\logand}{\wedge}

\title{Convergence rates and source conditions for Tikhonov regularization with sparsity constraints
  }
\author{Dirk A. Lorenz\thanks{Dirk Lorenz, Zentrum f\"ur Technomathematik, Fachbereich 3, Universit\"at Bremen, PO Box 330440, 28334 Bremen, Germany, \protect\url{dlorenz@math.uni-bremen.de}.}}

\begin{document}
\maketitle

\begin{abstract}
  This paper addresses the regularization by sparsity constraints by means of weighted $\ell^p$ penalties for $0\leq p\leq 2$.
  For $1\leq p\leq 2$ special attention is payed to convergence rates in norm and to source conditions.
  As main results it is proven that one gets a convergence rate of $\sqrt{\delta}$ in the 2-norm for $1< p\leq 2$ and in the 1-norm for $p=1$ as soon as the unknown solution is sparse.
  The case $p=1$ needs a special technique where not only Bregman distances but also a so-called Bregman-Taylor distance has to be employed.
  
  For $p<1$ only preliminary results are shown.
  These results indicate that, different from $p\geq 1$, the regularizing properties depend on the interplay of the operator and the basis of sparsity.
  A counterexample for $p=0$ shows that regularization need not to happen.
\end{abstract}


\section{Introduction}
\label{sec:introduction}
In this paper we discuss the regularizing properties of so-called sparsity constraints.
We consider linear inverse problems with a bounded operator $A:X\to Y$ between two Hilbert spaces.
Our setting is classical~\cite{engl1996inverseproblems}: We assume that we are given noisy data $g^\delta\in Y$ such that there exists
$g^+ = Af^+$ with $\norm{g^+-g^\delta}_Y\leq \delta$.
Our aim is to reconstruct $f^+$ from the noisy data $g^\delta$.
It is well known that this problem is ill-posed if and only if the range of $A$ is non-closed~\cite{engl1996inverseproblems}.

Recently regularization with sparsity constraints has become popular due to the influential paper~\cite{daubechies2003iteratethresh}.
In this setting one assumes, that the unknown solution has a sparse representation
in a certain orthonormal basis or frame $(\psi_k)$ of $X$,
i.e.~the unknown solution $f^+$ can be expressed as
$f^+ = \sum u_k \psi_k$ where the sum consists of a few
(and especially finitely many) terms only.
This knowledge is used to set up a so-called sparsity constraint for Tikhonov regularization,
i.e.~the regularized solution is given as a minimizer of
\begin{equation*}
  \norm{Af - g^\delta}_Y^2 + \alpha \sum_k w_k\phi(|\scp{f}{\psi_k}|)
\end{equation*}
with a suitably chosen function $\phi$.
The parameter $\alpha>0$ is a regularization parameter
and the weighting sequence $w_k>0$ allows to regularize each coefficient individually.
For the weighting sequence we assume that it is bounded away from zero: $w_k\geq w_0>0$.
Several choices of $\phi$ are possible.
In~\cite{daubechies2003iteratethresh} it is argued that the choice $\phi(s) = s^p$ for $1\leq p < 2$ promotes sparsity of the minimizer.
A heuristic explanation is that this functions give a higher weight to small coefficients and lower weight to large coefficients.
Of course the cases $p<1$ or even $p=0$ will produce sparse minimizers but in this case the convexity of the functional is lost and minimizers need not to exist (see~\cite{lorenz2007softhard} for a discussion of the case $A=I$).

For notational convenience we introduce the synthesis operator $B:\lpspace{2}\to X$ defined by $Bu = \sum_k u_k\psi_k$.
We define $K = AB$ and rewrite the Tikhonov functional as
\begin{equation}
  \label{eq:sparsity_funct}
  \Psi(u) = \norm{Ku-g^\delta}_Y^2 + \alpha \sum w_k |u_k|^p.
\end{equation}
The calculation of a minimizer of the above functional is not a straightforward task.
Convergent algorithms in the infinite dimensional setting for $1\leq p\leq 2$ were proposed and analyzed in \cite{daubechies2003iteratethresh,bredies2005gencondgrad,bredies2008harditer,daubechies2007projgrad,bredies2008itersoftconvlinear,griesse2008ssnsparsity}.
Generalizations to joint sparsity \cite{fornasier2008jointsparsity}, nonlinear operators \cite{ramlau2006tikhproject,ramlau2005tikhreplace,bonesky2007gencondgradnonlin} and the case $p=0$~\cite{blumensath2007iterhard} have been proposed.

In this paper we are going to discuss the regularizing properties of sparsity constraints.
First results on this topic can be found in~\cite{daubechies2003iteratethresh}
where convergence of the minimizers in $X$ (resp~$\ell^2$) for vanishing noise and the parameter choice $\alpha(\delta)$ such that $\alpha\to 0$ and $\delta^2/\alpha\to 0$ has been shown.
Moreover, it is shown that, in the special case of wavelet bases with a special class of weights which lead to Besov spaces, convergence rates can be achieved.
The paper~\cite{ramlau2006tikhproject} also deals with convergence of the minimizers and the proofs there show that convergence in the stronger $\ell^1$ norm holds.
Sparsity constraints can also be discussed in the framework of regularization in Banach spaces like, e.g., in~\cite{burger2004convarreg,resmerita2005regbanspaces,resmerita2006nonquadreg,hofmann2007convtikban,burger2007errorinvscalespace}.
In these papers convergence rates for general convex regularization are given in terms of Bregman distances.
In this paper we focus on convergence rates for sparsity constraints in norm, i.e.~in the norm in $X$ resp.~$\ell^2$ or the $\ell^1$-norm.

The paper is organized as follows.
Section~\ref{sec:preliminary-results} presents auxiliary results and in Section~\ref{sec:regul-with-1p2} results on convergence rates
for Tikhonov regularization with~(\ref{eq:sparsity_funct}) for $1<p\leq 2$ are presented,
especially we illustrate the role of the source condition.
Section~\ref{sec:regul-with-p=1} treats the case $p=1$ which is considerably different and a different technique has to be used.
The Section~\ref{sec:regul-with-p=0} collects preliminary results
on the regularization with $p<1$.
Here, no convergence rates can be given so far, and are not to be expected in general.
In the last section we draw conclusions.

\paragraph{Notation.}
\label{sec:notation}
We denote with $\ell^p_w$ the weighted $\ell^p$ space,
i.e.~the sequences $u$ such that $\sum w_k|u_k|^p$ converges.
We consider the spaces $\ell^p_w$ for $0< p\ < \infty$
which are normed spaces (quasi-normed for $p<1$)
when equipped with the (quasi-)norm
$\norm{u}_{p,w} = (\sum w_k|u_k|^p)^{1/p}$.
By $\ell^0$ we denote the set
$\{u:\N\to\R\ :\ u_k\not=0\ \text{ for finitely many }\ k\}$
of finitely supported or sparse sequences and with
$\ell^0_w$ the set
$\{u:\N\to\R\ :\ \sum w_k\sgn(|u_k|)<\infty \}$.
For simplicity we write $\norm{u} = \norm{u}_2$
and the inner product of $u,v\in\ell^2$ is denoted by $\scp{u}{v}$.
Moreover, we will frequently use component-wise application of operators to sequences, e.g. $(|u|^p)_k = |u_k|^p$ or $(wu)_k = w_k u_k$. 
With
\[
\Sgn(x) =
\begin{cases}
  \{1\} & \text{ for }\ x>0\\
  [{-1},1] & \text{ for }\ x=0\\
  \{{-1}\} & \text{ for }\ x<0.
\end{cases}
\]
we denote the multivalued sign while $\sgn$ stands for the usual sign with $\sgn(0)=0$.
For an operator $A:X\to Y$ between two Hilbert spaces the Hilbert space adjoint is denoted by $A^*:Y\to X$.

\section{Preliminary results}
\label{sec:preliminary-results}
In this section we collect preliminary results which are needed in the following.

As a first result we report that the cases $1\leq p<2$ indeed promote sparsity and that $p=1$ lead to finitely supported minimizers.
\begin{lemma}
  Let $1\leq p\leq 2$.
  A minimizer $u^*$ of $\Psi$ from~(\ref{eq:sparsity_funct}) fulfills
  \[u^* \in \ell^{2(p-1)}_{w^2}.\]
\end{lemma}
\begin{proof}
  Every minimizer $u$ of $\Phi$ fulfills
  \begin{equation}
    \label{eq:charac_minimizer_inclusion}
    -2K^*(Ku-g^\delta) \in \alpha w p \Sgn(u)|u|^{p-1}.
  \end{equation}
For $p>1$ the inclusion becomes an equation and since the left hand side is an $\ell^2$ sequence, the right hand side is also in $\ell^2$.
  It follows that
  \[\sum w_k^2 |u_k|^{2(p-1)} < \infty.\]
  For $p=1$ assume that $u\notin \ell^0_{w^2}$ i.e.~the sum $\sum w_k^2\sgn(|u_k|)$ diverges.
  Hence, every other choice of a sign in (\ref{eq:charac_minimizer_inclusion}) also leads to a diverging sum and it follows that the left hand side in~(\ref{eq:charac_minimizer_inclusion}) can not be an $\ell^2$ sequence, which is a contradiction.
\end{proof}

The next statement is on convergence of minimizers of~(\ref{eq:sparsity_funct}) for $\delta\to 0$.
\begin{theorem}[\cite{daubechies2003iteratethresh}]
  \label{thm:conv_minimizer}
  Assume that either $p>1$ or $K$ is injective, $w_k\geq w_0>0$, and let $u^{\alpha,\delta}$ be a minimizer of $\Psi$ from~(\ref{eq:sparsity_funct}). If the parameter choice $\alpha(\delta)$ fulfills
  \[
  \lim_{\delta\to 0} \alpha(\delta) = 0,\quad \lim_{\delta\to 0} \frac{\delta^2}{\alpha(\delta)} = 0
  \]
  then it holds
  \[
  \lim_{\delta\to 0}\norm{u^{\alpha,\delta} - u^+} =  0.
  \]
\end{theorem}
This says that that the method is indeed a regularization.
To get a statement on the rate of convergence the true solution $u^+$ has to fulfill some source condition.
This will be topic of sections~\ref{sec:regul-with-1p2} and~\ref{sec:regul-with-p=1}.

Next we state a basic inequality which we will need in the following.
\begin{lemma}[\cite{bredies2008harditer}]
  \label{lem:p_estimate}
  Let $1<p\leq 2$. For $C>0$ and $L>0$ it holds  for every
  $s,t\in\R$ with $|s|\leq C$ and $|t-s|\leq L$
  \[
  |t|^p - |s|^p\geq p\sgn(s)|s|^{p-1}(t-s) + \kappa|t-s|^2
  \]
  with $\kappa= \frac{p(p-1)}{2(C+L)^{2-p}}$.
\end{lemma}

\section{Regularization with $1<p\leq 2$}
\label{sec:regul-with-1p2}
In this section we analyze the ``easiest'' case $1<p\leq 2$.
The main result goes as follows.
\begin{theorem}
  \label{thm:convergence_rates}
  Let $1<p\leq 2$, $w_k\geq w_0>0$ and let $u^{\alpha,\delta}$ be a minimizer of $\Psi$ given in~(\ref{eq:sparsity_funct}).
   Furthermore let $u^+$ fulfill the source condition
  \begin{equation}\label{eq:source_condition}
    \exists \theta\in Y:\ w \sgn(u^+)|u^+|^{p-1} = K^*\theta.
  \end{equation}
  Then for the choice $\alpha \sim \delta$ it holds
  \begin{align}
    \label{eq:est_data-side}
    \norm{Ku^{\alpha,\delta} - g^\delta}_Y &= \bigO(\delta)\ \text{ for }\ \delta\to 0\\
    \label{eq:est_reconstruction-side}
    \sum w_k |u^{\alpha,\delta}_k - u^+_k|^2 &=  \bigO(\delta)\ \text{ for }\ \delta\to 0.    
  \end{align}
\end{theorem}
\begin{proof}
  Due to the minimizing property we have
  \[
  \norm{Ku^{\alpha,\delta}-g^\delta}_Y^2 +\alpha\sum w_k|u^{\alpha,\delta}_k|^p \leq 
    \delta^2 +\alpha\sum w_k|u^+_k|^p
  \]
  which gives
  \[
  \norm{Ku^{\alpha,\delta}-g^\delta}^2_Y +\alpha\sum w_k\bigl(|u^{\alpha,\delta}_k|^p - |u^+_k|^p\bigr) \leq 
  \delta^2.
  \]
  Since $|u^+_k|$ and $|u^{\alpha,\delta}_k-u^+_k|$ can be bounded uniformly in $k$ (the second due to Theorem~\ref{thm:conv_minimizer}) we can apply Lemma~\ref{lem:p_estimate} which yields
  \[
  \norm{Ku^{\alpha,\delta}-g^\delta}^2_Y +\alpha\kappa\sum w_k|u^{\alpha,\delta}_k - u^+_k|^2  + p\alpha\sum w_k\sgn(u^+_k)|u^+_k|^{p-1}(u^{\alpha,\delta}_k - u^+_k)\leq 
  \delta^2.  
  \]
  Rearranging gives
  \[
  \norm{Ku^{\alpha,\delta}-g^\delta}^2_Y + \alpha \kappa\sum w_k|u^{\alpha,\delta}_k - u^+_k|^2 \leq
    \delta^2 + \alpha\scp{p w \sgn(u^+)|u^+|^{p-1}}{u^+-u^{\alpha,\delta}}.
  \]
  Applying the source condition~(\ref{eq:source_condition}) and the Cauchy-Schwarz inequality leads to
  \[
  \norm{Ku^{\alpha,\delta}-g^\delta}^2_Y + \alpha \kappa\sum w_k|u^{\alpha,\delta}_k - u^+_k|^2 \leq
    \delta^2 + \alpha p\norm{\theta}_Y\norm{K(u^+-u^{\alpha,\delta})}_Y.
  \]
  Adding and subtracting $g^\delta$ in the last norm and denoting $\rho = \norm{\theta}_Y p/2$ leads to
  \[
  \norm{Ku^{\alpha,\delta}-g^\delta}^2_Y + \alpha \kappa\sum w_k|u^{\alpha,\delta}_k - u^+_k|^2 \leq
    \delta^2 + 2\alpha\rho\delta + 2\alpha\rho\norm{Ku^{\alpha,\delta} - g^\delta}_Y.
  \]
  Rearranging and completing the squares gives
  \[
  (\norm{Ku^{\alpha,\delta} - g^\delta}_Y - \alpha\rho)^2 + \alpha \kappa\sum w_k|u^{\alpha,\delta}_k - u^+_k|^2 \leq (\delta + \alpha\rho)^2.
  \]
  This finally implies
  \begin{equation}
    \label{eq:est_data_side_alpha}
  \norm{Ku^{\alpha,\delta}- g^\delta}_Y \leq \delta + 2\alpha\rho
  \end{equation}
  and
  \begin{equation}
    \label{eq:est_reconstruction_side_alpha}
  \sum w_k|u^{\alpha,\delta}_k - u^+_k|^2 \leq \frac{(\delta + \alpha\rho)^2}{\alpha \kappa}.
  \end{equation}
  The assertion follows with $\alpha \sim \delta$.
\end{proof}

Since $w_k\geq w_0$ we can deduce the following corollary immediately.
\begin{corollary}
  Under the assumptions of Theorem~\ref{thm:convergence_rates} it holds
  \[
  \norm{u^{\alpha,\delta}-u^+} = \bigO(\sqrt{\delta}).
  \]
\end{corollary}

We state a few remarks to illustrate Theorem~\ref{thm:convergence_rates}.

\begin{remark}[Constants in the $\bigO$-notation]
  \label{rem:constants_bigO}
  From~(\ref{eq:est_data_side_alpha}) one deduces that
  \[
  \norm{Ku^{\alpha,\delta}- g^\delta}_Y \leq (1+2\rho)\delta
  \]
  and hence, the constant in the $\bigO$ notation only depends on $\rho$.
  From the estimate~(\ref{eq:est_reconstruction_side_alpha}) we have
  \[
  \sum w_k|u^{\alpha,\delta}_k - u^+_k|^2 \leq \frac{(1 + \rho)^2}{\kappa}\delta.
  \]
  In this case, the constant depends also on $\kappa$ from Lemma~\ref{lem:p_estimate} for which it holds
  \[
  \frac{1}{\kappa} = \frac{2(C+L)^{2-p}}{p(p-1)}
  \]
  where $C$ is an upper bound on $|u^+_k|$ and $L$ is an upper bound for $|u^{\alpha,\delta}_k - u^+_k|$.
  The value $L$ tends to zero for $\delta\to 0$ and $C$ depends on $u^+$ only and hence, $C$ and $L$ are uniformly bounded for $\delta\to 0$.
  Finally, we see that the constant $1/\kappa$ mainly depends on $p$
  and is large for small $p$
  and namely it tends to infinity for $p\to 1$.
  To summarize, we may say that the regularization with a weighted $\ell^p$-norm leads to a convergence rate of order $\sqrt{\delta}$ in the 2-norm but the associated constant gets arbitrarily large for $p$ close to one.
  Hence, one may not assume a similar theorem to hold for the limiting case $p=1$.
  Fortunately, Theorem 4.3 below shows that this pessimism unfounded.
\end{remark}

\begin{remark}[The results of Burger and Osher~\cite{burger2004convarreg}]
  \label{rem:burger_osher}
  In the case of a general convex and lower-semicontinuous penalty functional $J$, Burger and Osher proved that the source condition
  \begin{equation*}
    \exists \theta:\ K^*\theta \in \partial J(u^+)
  \end{equation*}
  leads to a convergence rate
  \[
  D_\xi(u^{\alpha,\delta},u^+) = \bigO(\delta)
  \]
  for $u^{\alpha,\delta}$ minimizers of
  \[
  \norm{Ku-g^\delta}_Y + \alpha J(u).
  \]
  Here $\partial J$ denotes the subgradient of $J$, $\xi\in\partial J(u^+)$
  and 
  \[D_\xi(u^{\alpha,\delta},u^+) = J(u^{\alpha,\delta}) - J(u^+) - \scp{\xi}{u^{\alpha,\delta} - u^+}
  \]
  is the Bregman distance, see also \cite{hofmann2007convtikban}.
  One can also deduce Theorem~\ref{thm:convergence_rates} from this result by noting that this source condition is precisely the one in Theorem~\ref{thm:convergence_rates} and that for $1<p\leq 2$ the Bregman distance of $J(u) = \sum w_k|u_k|^p$ can be bounded from below:
  \[
  \sum w_k|u^{\alpha,\delta}-u^+|^2 \leq D_\xi(u^{\alpha,\delta},u^+)
  \]
  for $\norm{u^{\alpha,\delta}-u^+}<M$
  (which follows from Lemma~\ref{lem:p_estimate} or the inequalities of Xu and Roach, see~\cite{xu1991charinequalities,schoepfer2006illposedbanach,schuster2008tikbanach}).
\end{remark}

\begin{remark}[Source conditions in terms of $\lpspace{p}$-spaces]
  \label{rem:source_cond_lp}
  In the classical (quadratic) theory the source condition can usually be interpreted as some kind of smoothness condition.
  When working in sequence space, we see that the source condition~(\ref{eq:source_condition}) says something about the decay of the solution $u^+$.
  We assume that the operator under consideration has  the property
  $\range K^* = \lpspace{q}_v$ where we assume that the space
  $\lpspace{q}_v$
  is contained in $\lpspace{2}$.  
  Hence, the dual space $(\lpspace{q}_v)' = \lpspace{q'}_{v'}$
  with dual exponent $q' = q/(q-1)$ and dual weight $v^{-1/(q-1)}$
  is larger than $\lpspace{2}$.
  One may say that the operator
  $K:\lpspace{q'}_{v'}\to Y$ has a ``smoothing''
  (or better ``damping'') property.
  Now, the source condition~(\ref{eq:source_condition}) reads as
  $w\sgn(u^+)|u^+|^{p-1} \in \range K^* = \lpspace{q}_v$ and hence
  \[
  \sum v_k w_k^q |u^+_k|^{q(p-1)} < \infty\ \text{ or equivalently }\
  u^+ \in \lpspace{q(p-1)}_{vw^q}.
  \]
\end{remark}
\section{Regularization with $p=1$}
\label{sec:regul-with-p=1}
We now turn to the case $p=1$.
In this case previous results give convergence rates in the Bregman distance only~\cite{burger2004convarreg,resmerita2005regbanspaces,resmerita2006nonquadreg,hofmann2007convtikban}.
Moreover, Remark~\ref{rem:burger_osher} does not apply, since the function
$J(u) = \sum_k w_k|u_k|$ is not strictly convex
and hence, the Bregman distance with respect to the functional
$J(u) = \sum w_k|u_k|$ can not be estimated
by the $\lpspace{2}$-norm in general.
It holds $\partial J(u) = (w_k\Sgn(u_k))_k$.
One sees that the Bregman distance fulfills
\[
D_\xi(u,u^+) \leq 2 \sum_{(u_k>0 \logand u_k\leq 0)\atop \logor (u_k<0 \logand u_k\geq 0)}|u_k|.\]
Consequently, the Bregman distance is zero as soon as the signs of $u$ and $u^+$ coincide and a convergence rate regarding the Bregman distance does not give satisfactory information, see also~\cite{burger2007errorinvscalespace}.

To prove a convergence rate like in Theorem~\ref{thm:convergence_rates} we need the following lemma which can be found in similar form in~\cite{bredies2008itersoftconvlinear}.
As an important ingredient we need the so called FBI property, also from~\cite{bredies2008itersoftconvlinear}.
\begin{definition}
  \label{def:fbi}
  An operator $K:\lpspace{2}\to Y$ mapping into a Hilbert
  space has the
  \emph{finite basis injectivity} (FBI) property, if for all finite subsets
  $I\subset\N$ the operator $K|_I$ is injective, i.e.~for
  all $u,v\in \ell^2$ with $Ku = Kv$ and
  $u_k = v_k = 0$ for all $k \notin I$ it follows $u = v$.
\end{definition}
The lemma gives an estimate which compares the Bregman distance with the $\ell^1$-norm.
\begin{lemma}
  \label{thm:bregman_taylor_dist}
  Let $u^+$ have finite support, $w_k\geq w_0>0$, let $K$ fulfill the FBI property, and define
  \begin{eqnarray}
    \label{eq:def_taylor_distance}
    T(u) & = & \norm{K(u-u^+)}_Y^2\\
    \label{eq:def_bregman_distance}
    R(u) & = & \sum w_k|u_k| - \sum w_k|u^+_k| - \sum w_k\sgn(u^+_k)(u_k-u^+_k).
  \end{eqnarray}
  Then there exists $\lambda>0$ such that
  \[
  R(u) + T(u) \geq \lambda \norm{u-u^+}_1^2
  \]
  whenever $\norm{u-u^+}_1\leq M$.
\end{lemma}
\begin{proof}
  We define
  $I = \{k\ |\ \sgn(u^+_k) = \pm 1\}$
  which is a finite set.
  We estimate
  \begin{align*}
    R(u) & =  \sum_k w_k|u_k| - w_k|u^+_k| - w_k\sgn(u^+_k)(u_k-u^+_k)\\
    & = \sum_k w_k|u_k| - w_k\sgn(u^+_k)u_k\\
    & \geq \sum_{k\notin I} w_k|u_k| - w_k\sgn(u^+_k)u_k = \sum_{k\notin I}w_k |u_k|.
  \end{align*}
  Denoting with $I^c$ the complement of $I$ and with $P_{I^c}$ the projection onto the subspace where all coefficients in $I$ are zero we get (using $u^+_k=0$ for $k\in I^c$)
  \[
  R(u) \geq w_0\norm{P_{I^c}(u - u^+)}_1.
  \]
  Since $\norm{P_{I^c}(u-u^+)} \leq M $ we can estimate
  \begin{equation}
    \label{eq:Bregman_dist_from_below}
    R(u) \geq \frac{w_0}{M}\norm{P_{I^c}(u - u^+)}_1^2.
  \end{equation}
  
  To establish an estimate for the remaining part $P_Iu$ we start
  with $u = P_{I}u + P_{I^c}u$ and use the inequalities of Cauchy-Schwarz
  (in the form $-\scp{u}{v} \leq \norm{u}\norm{v}$)
  and Young ($ab\leq \tfrac{a^2}{4}+b^2$ for $a,b>0$) to get
  \begin{align}
    \norm{Ku}_Y^2 & = \norm{KP_{I}u}_Y^2 + 2\scp{KP_{I}u}{KP_{I^c} u} + \norm{KP_{I^c} u}_Y^2\nonumber\\
    & \geq \frac{\norm{KP_I u}_Y^2}{2} - \norm{KP_{I^c}u}_Y^2\nonumber\\
    & \geq \frac{\norm{KP_I u}_Y^2}{2} - \norm{K}^2\norm{P_{I^c}u}^2.\label{eq:est_Ku_below}
  \end{align}
  Since $I$ is finite and $K$ obeys the FBI property
  there is a constant $c>0$ such that
  \[
  c\norm{P_I u}^2 \leq \norm{KP_I u}_Y^2.
  \]
  Moreover, again since $I$ is finite, we can estimate the 2-norm from below by the 1-norm which leads to
  \[
  \tilde c\norm{P_I u}_1^2 \leq \norm{KP_I u}_Y^2.
  \]
  Combining this with~(\ref{eq:est_Ku_below}) gives
  \[
  \norm{P_I u}_1^2 \leq \frac{2}{\tilde c}(\norm{Ku}_Y^2 + \norm{K}^2\norm{P_{I^c}u}^2) 
  \]
  Applying this estimate to $u-u^+$ instead of $u$ and adding the inequality~(\ref{eq:Bregman_dist_from_below}) leads to
  \[
  \norm{u-u^+}_1^2 \leq \frac{2}{\tilde c}(T(u) + \norm{K}^2\norm{P_{I^c}(u-u^+)}^2) + \frac{M}{w_0}R(u).
  \]
  By estimating the 1-norm from below by the 2-norm in~(\ref{eq:Bregman_dist_from_below}) we get $\tfrac{M}{w_0}R(u) \geq \norm{P_{I^c}(u - u^+)}_2^2$ and hence,
  \[
  \norm{u-u^+}_1^2 \leq \frac{2}{\tilde c} T(u) + \frac{M}{w_0}\bigl(\frac{2\norm{K}^2}{\tilde c} + 1\bigr)R(u)
  \]
  which proves the claim.
\end{proof}

While the term $R$ from~(\ref{eq:def_bregman_distance}) is a Bregman distance, the term $T$ from~(\ref{eq:def_taylor_distance}) can be seen as Taylor distance: We define the functional $F(u) = \norm{Ku-g^\delta}^2_Y$ and observe that the term $T$ can be rewritten as
\[
T(u) = F(u) - F(u^+) - \scp{F'(u^+)}{u-u^+}.
\]
Consequently, $T$ is the remainder of the Taylor expansion of the fidelity term $F$.
Therefore, Lemma~\ref{thm:bregman_taylor_dist} can be seen as an estimate on the Bregman-Taylor-distance $R+T$.

Lemma~\ref{thm:bregman_taylor_dist} enables us to prove the main result of this paper:

\begin{theorem}
  \label{thm:reg_p1}
  Let $u^+$ have finite support, $w_k\geq w_0>0$, $K$ obey the FBI property, and let furthermore
  $u^+$ fulfill the source condition
  \begin{equation}
    \label{eq:source_cond_p1}
    \exists \theta\in Y:\ w\sgn(u^+)= K^*\theta.
  \end{equation}
  Then for every
  \[
  u^{\alpha,\delta} \in \argmin\ \norm{Ku-g^\delta}^2_Y + \alpha\sum w_k|u_k|
  \]
  it holds
  \[
  \norm{u^{\alpha,\delta}-u^+}_1 = \bigO(\sqrt{\delta}).
  \]
\end{theorem}
\begin{proof}
  Due to the minimizing property we have
  \begin{align*}
    0 & \leq \norm{Ku^+ - g^\delta}^2_Y + \alpha\sum_k w_k|u^+_k|
         - \norm{Ku^{\alpha,\delta} - g^\delta}^2_Y - \alpha\sum_k w_k|u^{\alpha,\delta}_k|\\
      & = \norm{Ku^+ - g^\delta}^2_Y - \norm{Ku^{\alpha,\delta} - g^\delta}^2_Y\\
      &\quad {}+\alpha( \sum_k w_k|u^+_k| - \sum_k w_k|u^{\alpha,\delta}_k| + \sum_k w_k\sgn(u^+_k)(u^{\alpha,\delta}_k - u^+_k))\\
      &\quad {}- \alpha\sum_k w_k\sgn(u^+_k)(u^{\alpha,\delta}_k - u^+_k).
  \end{align*}
  Rearranging gives
  \begin{equation*}
    \alpha R(u^{\alpha,\delta}) \leq \delta^2 - \norm{Ku^{\alpha,\delta} - g^\delta}^2_Y- \alpha\sum_k w_k\sgn(u^+_k)(u^{\alpha,\delta}_k - u^+_k).
  \end{equation*}
  Since the convergence $u^{\alpha,\delta}\to u^+$ is known from Theorem~\ref{thm:conv_minimizer} we can use Lemma~\ref{thm:bregman_taylor_dist} to obtain
  \[
  \alpha\lambda\norm{u^{\alpha,\delta}-u^+}_1^2 - \alpha\norm{K(u^{\alpha,\delta}-u^+)}^2_Y \leq
  \delta^2 - \norm{Ku^{\alpha,\delta} - g^\delta}^2_Y- \alpha\sum_k w_k\sgn(u^+_k)(u^{\alpha,\delta}_k - u^+_k).
  \]
  With the source condition~(\ref{eq:source_cond_p1}),
  the notation $\rho = \norm{\theta}_Y/2$,
  and the Cauchy-Schwarz inequality this gives
  \[
  \alpha\lambda\norm{u^{\alpha,\delta}-u^+}_1^2 - \alpha\norm{K(u^{\alpha,\delta}-u^+)}^2_Y \leq
  \delta^2 - \norm{Ku^{\alpha,\delta} - g^\delta}^2 + \alpha 2\rho\norm{K(u^{\alpha,\delta} - u^+)}.
  \]
  Adding and subtracting $g^\delta$ in the last norm and rearranging leads to
  \[
  \alpha\lambda\norm{u^{\alpha,\delta}-u^+}_1^2 - \alpha\norm{K(u^{\alpha,\delta}-u^+)}^2_Y  + \norm{Ku^{\alpha,\delta} - g^\delta}^2_Y - 2\alpha \rho\norm{Ku^{\alpha,\delta} - g^\delta}_Y\leq
  \delta^2 + 2\alpha \rho\delta.
  \]
  Using
  \[
  \norm{K(u^{\alpha,\delta}-u^+)}^2_Y \leq \norm{Ku^{\alpha,\delta}-g^\delta}^2_Y + 2\delta\norm{Ku^{\alpha,\delta} - g^\delta}_Y +\delta^2
  \]
  leads to
  \[
  \alpha\lambda\norm{u^{\alpha,\delta}-u^+}_1^2 + (1-\alpha)\norm{K(u^{\alpha,\delta}-g^\delta)}^2_Y - 2\alpha(\rho+\delta)\norm{Ku^{\alpha,\delta}-g^\delta}_Y \leq (1+\alpha)\delta^2 + 2\alpha\rho\delta.
  \]
  Dividing by $(1-\alpha)$ and completing the square on the left hand side gives
  \[
  \frac{\alpha}{1-\alpha}\lambda\norm{u^{\alpha,\delta}-u^+}_1^2
    + \Bigl(\norm{Ku^{\alpha,\delta}-g^\delta}_Y - \frac{\alpha}{1-\alpha}(\rho+\delta)\Bigr)^2
    \leq
    \frac{1+\alpha}{1-\alpha}\delta^2 + \frac{2\alpha\rho\delta}{1-\alpha} + \Bigl(\frac{\alpha}{1-\alpha}\Bigr)^2(\rho+\delta)^2.
  \]
  Finally, this gives
  \begin{align}
  \norm{u^{\alpha,\delta}-u^+}_1^2 & \leq \frac{1}{\lambda}\Bigl( \frac{1+\alpha}{\alpha}\delta^2 + 2\rho\delta + \frac{\alpha}{1-\alpha}(\rho+\delta)^2 \Bigr)\nonumber\\
  & = \frac{1}{\lambda\alpha(1-\alpha)}\Bigl( \delta + \alpha\rho \Bigr)^2.
  \label{eq:explicit_p1}
  \end{align}
  The choice $\alpha = \delta$ proves
  \[
  \norm{u^{\alpha,\delta}-u^+}_1^2 = \bigO(\delta)\ \text{ for }\ \delta\to 0.
  \]

\end{proof}

For $p=1$ the source condition says that $u^+$ must only have a finite number of non-zero entries.
This is the natural limit for $p\to 1$ as can be seen from Remark~\ref{rem:source_cond_lp}.

Theorem~\ref{thm:reg_p1} is remarkable since, as mentioned in Remark~\ref{rem:constants_bigO}, the constant in the $\bigO$-notation in Theorem~\ref{thm:convergence_rates} blows up to infinity for $p\to 1$.
Equation~(\ref{eq:explicit_p1}) shows that the constant in the $\bigO$-notation depends on the constant $\lambda$ from Lemma~\ref{thm:bregman_taylor_dist} and on $\rho = \norm{\theta}_Y/2$ only.
Basically the constant $1/\kappa$ in Remark~\ref{rem:constants_bigO} has been replaced by $1/\lambda$ from Lemma~\ref{thm:bregman_taylor_dist}.

\begin{remark}[The result of Hofmann et al.~\cite{hofmann2007convtikban}]
Hofmann et al.~considered in~\cite{hofmann2007convtikban} general convex regularization of operator equations in Banach spaces of the form
\[
\norm{F(u)-g^\delta}_Y^p + \alpha J(u).
\]
They showed a convergence rate of $\bigO(\delta)$ in the Bregman distance for non-smooth operators $F$ under the source condition that there exists $\beta_1\in[0,1[$, $\beta_2\geq 0$ and $\xi\in\partial J(u^+)$ such that
\[
-\scp{\xi}{u-u^+} \leq \beta_1 D_\xi(u,u^+) + \beta_2\norm{F(u)-F(u^+)}
\]
(note that the negative sign on the left hand side is a typo in the original paper).
This source condition is difficult to check in concrete situations.
Applied to the situation of Theorem~\ref{thm:reg_p1} it reads as:
There exists $\xi\in w\Sgn(u^+)$ such that
\[
-\scp{\xi}{u-u^+} \leq \beta_1 D_\xi(u,u^+) + \beta_2\norm{K(u-u^+)}_Y.
\]
This condition is for example fulfilled if the sequence $w_k$ is bounded and
\[
\norm{u-u^+}_1 \leq\frac{1}{\max w_k} (\beta_1 R(u) + \beta_2\norm{K(u-u^+)}_Y)
\]
which resembles the Bregman-Taylor estimate from Lemma 4.2.
However, Theorem 4.3 gives a convergence rate of $\bigO(\sqrt{\delta})$ in the $\ell^1$-norm and the Bregman-Taylor estimate is only needed to pass from the Bregman distance to the $\ell^1$-norm.
Additionally, Theorem 4.3 needs the source condition~(\ref{eq:source_cond_p1}).
  
\end{remark}

\section{Regularization with $p<1$?}
\label{sec:regul-with-p=0}
The functional~(\ref{eq:sparsity_funct}) is not convex if $p<1$.
Hence, there is no guarantee for uniqueness or existence of a minimizer.
In this section we show two extreme examples: One in which there exist minimizers which can be computed explicitly and regularization can be proven and the other where no minimizer exists at all.

\subsection{Regularization is possible}
\label{sec:regul-poss}
In this example we use an orthonormal basis which is perfectly adapted to the operator: the singular basis.
The singular value decomposition $(\sigma_k, \psi_k, \phi_k)$ of the operator $A$ consists of the singular values $\sigma_k$ and two orthonormal bases $\psi_k$ and $\phi_k$ of $X$ resp.~$Y$.
The operator $A$ can now be expressed as
\[Af = \sum_k \sigma_k \scp{f}{\psi_k}\phi_k.\]

Now we seek for a solution of $Af = g$ which is sparse in the basis $\psi_k$, i.e.~we have $u_k = \scp{f}{\psi_k}$ in (\ref{eq:sparsity_funct}).
Hence, the operator $K=AB$ has the form
\begin{equation}
  \label{eq:K_in_sing_basis}
Ku = ABu = \sum_k u_k\sigma_k \phi_k
\end{equation}

\begin{figure}[tb]
  \centering
  \includegraphics[width=.5\textwidth]{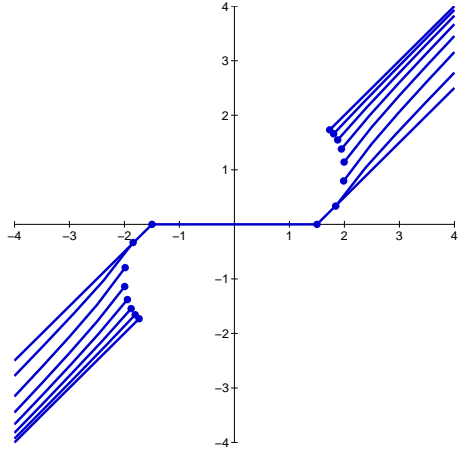}
  \caption{The thresholding functions $H^p_\alpha$ for $p=0, 0.15, 0.3, 0.45, 0.6, 0.75, 0.9, 1$ and $\alpha = 3$.}
  \label{fig:hardsoftinterp}
  
\end{figure}

To express the minimizer of~(\ref{eq:sparsity_funct}) we need the following function:
\begin{equation}
  \label{eq:Hp_def}
  H^p_\alpha(x) = \argmin_{y}\ (y-x)^2 + \alpha|y|^p.
\end{equation}
Note that this function can be multivalued in general.
The next lemma from~\cite{lorenz2007softhard} gives an implicit representation of the function $H^p_\alpha$.
\begin{lemma}
  \label{lem:Hp}
  Let 
  \begin{equation}
    \label{eq:def_G}
    G^p_\alpha(y) = y + \tfrac{\alpha p}{2} \sgn(y)|y|^{p-1}.
  \end{equation}
  The mapping $H^p_\alpha$ is given by the following formulae:
\begin{enumerate}
\item Let $1<p\leq 2$. Then $(G^p_\alpha)^{-1}$ exists and is single valued
      and it holds
      \[
      H^p_\alpha(x) = (G^p_\alpha)^{-1}(x).
      \]
\item Let $p=1$. Then
      \[H^1_\alpha(x) = \max(|x|-\alpha/2,0)\sgn(x).\]
\item Let $0< p<1$. Then
      \begin{equation}
        \label{eq:Hp_formula}
      H^p_\alpha(x) = \begin{cases}
                        0 & , \text{ for } |x| \leq \alphaeff\\
                        \parbox{.6\textwidth}
                            {the value of largest absolute value of the
                            inverse mapping of $G^p_\alpha$} & , \text{ for } |x|\geq
                            \alphaeff
                       \end{cases}
      \end{equation}
      where $\alphaeff = \frac{2-p}{2-2p}\Bigl(
      \alpha(1-p)\Bigr)^\frac{1}{2-p}$.
\item Let $p=0$. Then
  \begin{equation}
    \label{eq:eq:H0_formula}
    H^0_\alpha(x) =  
      \begin{cases}
        0 & , \text{ for } |x|\leq \alphaeff\\
        x & , \text{ for } |x|\geq \alphaeff\\
      \end{cases}
  \end{equation}
  where $\alphaeff = \sqrt{\alpha}$.
\end{enumerate}
\end{lemma}
The description in 3. may be a little unfamiliar.
For $p<1$ the function $G^p_\alpha$ is multivalued in $0$ with $G^p_\alpha(0) = \R$.
Its inverse is again multivalued (in fact it has at most three values) and the function $H^p_\alpha$ chooses either 0 or the value of largest absolute value, see Figure~\ref{fig:hardsoftinterp} and \cite{lorenz2007softhard} for more details.
Note moreover that for $p<1$ the function $H^p_\alpha$ is multivalued itself, namely it has two values for $|x| = \alpha_\text{eff}$.
For convenience we always choose the value 0 at these points in the following.

The next lemma is an easy consequence of the above lemma and the fact that the operator $K$ is diagonal with respect to the basis $(\phi_k)$ of $Y$.

\begin{lemma}
  Let $(\phi_k)$ be an orthonormal basis of $Y$ and let the operator $K:\lpspace{2}\to Y$ be given by~(\ref{eq:K_in_sing_basis}).
  Then, a minimizer of~(\ref{eq:sparsity_funct}) is given by
  \begin{equation}
    \label{eq:minimizer_sing_basis}
    u^{\alpha,\delta}_k  = 
    \begin{cases}
      \frac{1}{\sigma_k}H^p_{\alpha/\sigma_k^p}(\scp{g^\delta}{\phi_k}) & ,\text{ for } \sigma_k>0\\
      0 & ,\text{ for } \sigma_k=0.
      
    \end{cases} 
  \end{equation}
\end{lemma}

\begin{definition}
  For $0\leq p \leq 2$ we define the operator $R^p_\alpha: Y \to \lpspace{2}$ by
  \[
  \bigl(R^p_\alpha(g)\bigr)_k =
  \begin{cases}
    \tfrac{1}{\sigma_k}H^p_{\alpha/\sigma_k^p}(\scp{g^\delta}{\phi_k}) & ,\text{ for } \sigma_k>0\\
      0 & ,\text{ for } \sigma_k=0.
  \end{cases}
  \]
\end{definition}
Note that $R^p_\alpha$ is non-linear and discontinuous.

\begin{theorem}
  \label{thm:reg_p_leq_1}
  Let $0\leq p <1$. The operator $R^p_\alpha$ is
  \begin{enumerate}
  \item defined for every $g\in Y$.
  \item a regularization, i.e.~for $g\in\dom(K^+)$ it holds
    \[
    \lim_{\alpha\to 0} \norm{R^p_\alpha g - K^+ g} = 0.
    \]
  \end{enumerate}
\end{theorem}
\begin{proof}
  We abbreviate $g_k = \scp{g}{\phi_k}$.
  The pseudo-inverse is given by
  \[
  (K^+g)_k = 
  \begin{cases}
    g_k/\sigma_k &, \text{ for } \sigma_k>0\\
    0            &, \text{ for } \sigma_k=0
  \end{cases}
  \]
  and by the Picard condition
  this is an $\lpspace{2}$ sequence.
  For an $M\in\N$ we write
  \begin{align*}
    \norm{R^p_\alpha(g) - K^+g}^2 & = \sum_{\sigma_k>0} \frac{|H^p_{\alpha/\sigma_k^p}(g_k) - g_k|^2}{\sigma_k^2}\\
    & = \sum_{\sigma_k>0,\ k\leq M} \frac{|H^p_{\alpha/\sigma_k^p}(g_k) - g_k|^2}{\sigma_k^2} + 
    \sum_{\sigma_k>0,\ k> M} \frac{|H^p_{\alpha/\sigma_k^p}(g_k) - g_k|^2}{\sigma_k^2}.
  \end{align*} 
  For a given $\epsilon>0$ we choose $M$ such that
  $\sum_{\sigma_k>0,\ k> M} |g_k|^2/\sigma_k^2 < \epsilon$.
  Since we can deduce from Lemma~\ref{lem:Hp}
  \begin{equation*}
    |H^p_\alpha(x) - x| \leq |x|\label{eq:est_Hp_diff_pointwise}
  \end{equation*}
  we can estimate
  \begin{equation*}
    \norm{R^p_\alpha(g) - K^+g}^2  = \sum_{\sigma_k>0,\ k\leq M} \frac{|H^p_{\alpha/\sigma_k^p}(g_k) - g_k|^2}{\sigma_k^2} + 
    \epsilon.
  \end{equation*} 
  Furthermore, we see from Lemma~\ref{lem:Hp}
  \begin{equation*}
    H^p_\alpha(x) \to x\ \text{ for } \alpha\to 0.\label{eq:Hp_alpha_to_0}
  \end{equation*}
  and hence, for sufficiently small $\alpha$ we have
  \[
  \norm{R^p_\alpha(g) - K^+g}^2 < 2\epsilon.
  \]
\end{proof}

The above theorem does only proof convergence on the range of the operator.
To obtain results on the speed of convergence one may 
assume special sparseness or decay properties similar to~\cite{borup2003shrinkops}.
We are not going to pursue further in this direction
since the case of the singular basis is of limited interest
in practical applications.
Moreover, convergence for noisy data has not been shown.

\subsection{Regularization is impossible}
\label{sec:regul-imposs}
In this section we present an example where a sparsity constraint with exponent $p=0$ does not lead to a regularization.
In particular the minimization of the Tikhonov functional is not well-posed in the sense that it does not have a solution.
To this end, we design an operator $A$ which does not act well on a given
orthonormal basis $(\psi_k)$.
Let $\{h_k\}$ be a countable set which is dense in the unit-ball of $Y$,
i.e.~$\norm{h_k}_Y = 1$ and for every $g\in Y$ with $\norm{g}_Y=1$
and every $\epsilon>0$ there is an index $k_0$
such that $\norm{g-h_{k_0}}_Y\leq \epsilon$.
We define the operator $A$ on the basis $(\psi_k)$ by
\begin{equation}
  \label{eq:A}
  A\psi_k = h_k, \text{ i.e. } Ku = \sum_k u_k h_k.
\end{equation}
\begin{proposition}
\label{prop:reg_impossible}
Let $K$ be defined by~(\ref{eq:A}), $\norm{g}_Y^2> \alpha$ and let further $g$ be not a multiple of $h_k$ for every $k$.
Then the functional
\[\Psi(u) = \norm{Ku-g}^2_Y + \alpha\sum_k \sgn(|u_k|)\]
does not have a minimizer.
\end{proposition}
\begin{proof}
Since the penalty term $\sum_k \sgn(|u_k|)$ does only depend on the number of coefficients we minimize separately over subspaces of a given dimension $n$.

As first case we consider $n=0$, i.e.~we minimize just over $u=0$.
We observe that $\Psi(0) = \norm{g}^2_Y$.

As second case we observe that $\Psi(u)\geq 2\alpha$
if $u$ has more than two different non-zero entries.

The last case is to minimize over the one-dimensional subspaces
$X_k = \Span\{e_k\}$ where $e_k$ is the canonical basis of $\lpspace{2}$.
The values of $\Psi$ are
\[\Psi(d_k e_k) = \norm{d_k h_k - g}^2_Y + \alpha.\]
Since $\{h^k\}$ is dense in the unit ball may take $d_k = \norm{g}_Y$
and find a sequence $h_l$ such that $\norm{g}_Y h_l\to g$ for $l\to\infty$.
Hence, the minimal value of $\Psi$ over all subspaces $X_k$ is $\alpha$, i.e.
\[\inf_{u\in\bigcup X_k}\Psi(u) = \alpha\]
and this infimum is not attained since $g$ is not a multiple of a basis vector $h_k$.
\end{proof}
It is clear that a similar example can be constructed if the vectors $h_k$ accumulate at a single point: take $g$ as the accumulation point of $h_k$.

\begin{remark}
  We remark that also the constrained model
  \begin{equation}
    \label{eq:constrained}
    \tag{$P_\epsilon$}
    \text{Minimize } \sum_k\sgn(|u_k|) \text{ s.t. } \norm{Ku-g^\delta}_Y\leq \epsilon
  \end{equation}
  is not well posed with $K$ from~(\ref{eq:A})
  since it has an infinite number of solutions.
  One may say that this situation is a little better than that of
  Proposition~\ref{prop:reg_impossible} since now solutions are available.
  An easy example shows, that regularization need not to happen in this setup.
  Let $g^+ = h_1$ and let $\norm{g^+-g^\delta}_Y\leq\delta$.
  The corresponding true solution is $u^+=e_1$.
  Then there is a sequence $h_l$ such that $h_l\to h_1=g^+$.
  Moreover, for sufficiently large $l$, $u^{\epsilon,\delta} = e_l$ is a solution of~(\ref{eq:constrained}) with $\epsilon = \tau\delta$ with $\tau>1$
  (assumed that the norm of $g^\delta$ is not too small).
  Finally, $\norm{u^{\epsilon,\delta} - u^+} = \sqrt{2}$ is not converging to zero for $\epsilon = \tau\delta$ and $\delta\to 0$.
\end{remark}

\section{Conclusions}
\label{sec:conclusions}
In this paper the regularizing properties of sparsity constraints have been analyzed.
Special attention was payed to convergence rates in norm and to the source conditions.
For $1<p\leq 2$ we could show, as a simple application of the results of Burger, Osher~\cite{burger2004convarreg} and the inequality of Xu and Roach~\cite{xu1991charinequalities} (or the basic inequality in Lemma~\ref{lem:p_estimate} from~\cite{bredies2008harditer}), that a convergence rate $\sqrt{\delta}$ in the 2-norm can be achieved by a source condition saying that $u^+$ has to be in a weighted $\ell^p$ space with small $p$, see Remark~\ref{rem:source_cond_lp}.

The case $p=1$ needed a special technique: the Bregman-Taylor-distance from~\cite{bredies2008itersoftconvlinear}.
Applying this, a convergence rate $\sqrt{\delta}$ in the stronger 1-norm could be achieved under the source condition that $u^+$ is finitely supported.

The incipient discussion on regularization with $p<1$ showed two things:
First, regularization may or may not be possible and second, the regularization properties depend on the interplay of the operator $A$ and on the choice of the basis functions $(\psi_k)$---a phenomenon which is not known for $p\geq 1$.
One may conjecture that if the operator $A$ acts well on the basis $(\psi_k)$
(in the sense that the values $\frac{\scp{A\psi_k}{A\psi_j}}{\norm{A\psi_k}_Y\norm{A\psi_j}_Y}$ are not too large) regularization is possible.
This would parallel observations in the framework of compressed sensing on the mutual coherence of dictionaries, see~\cite{donoho2003sparse}.

\bibliographystyle{plain}
\bibliography{literature}

\begin{thebibliography}{10}

\bibitem{blumensath2007iterhard}
Thomas Blumensath, Mehrdad Yaghoobi, and Mike Davies.
\newblock Iterative hard thresholding and $l^0$ regularisation.
\newblock In {\em IEEE International Conference on Acoustics, Speech and Signal
  Processing}, April 2007.

\bibitem{bonesky2007gencondgradnonlin}
Thomas Bonesky, Kristian Bredies, Dirk~A. Lorenz, and Peter Maass.
\newblock A generalized conditional gradient method for nonlinear operator
  equations with sparsity constraints.
\newblock {\em Inverse Problems}, 23:2041--2058, 2007.

\bibitem{borup2003shrinkops}
Lasse Borup and Morten Nielsen.
\newblock Some remarks on shrinkage operators.
\newblock Aalborg University, Dept. of Mathematical Sciences, 2003.

\bibitem{bredies2008harditer}
Kristian Bredies and Dirk~A. Lorenz.
\newblock Iterated hard shrinkage for minimization problems with sparsity
  constraints.
\newblock {\em SIAM Journal on Scientific Computing}, 30(2):657--683, 2008.

\bibitem{bredies2008itersoftconvlinear}
Kristian Bredies and Dirk~A. Lorenz.
\newblock Linear convergence of iterated soft-thresholding.
\newblock {\em To appear in Journal of Fourier Analysis and Applications},
  2008.

\bibitem{bredies2005gencondgrad}
Kristian Bredies, Dirk~A. Lorenz, and Peter Maass.
\newblock A generalized conditional gradient method and its connection to an
  iterative shrinkage method.
\newblock To appear in \emph{Computational Optimization and Applications},
  2008.

\bibitem{burger2004convarreg}
Martin Burger and Stanley Osher.
\newblock Convergence rates of convex variational regularization.
\newblock {\em Inverse Problems}, 20(5):1411--1420, 2004.

\bibitem{burger2007errorinvscalespace}
Matrin Burger, Elena Resmerita, and Lin He.
\newblock Error estimation for {B}regman iterations and inverse scale space
  methods in image restoration.
\newblock {\em Computing}, 81(2--3):109--135, 2007.
\newblock Special Issue on Industrial Geometry (Guest editors: B. J\"uttler, H.
  Pottmann, O. Scherzer).

\bibitem{daubechies2003iteratethresh}
Ingrid Daubechies, Michel Defrise, and Christine De~Mol.
\newblock An iterative thresholding algorithm for linear inverse problems with
  a sparsity constraint.
\newblock {\em Communications in Pure and Applied Mathematics},
  57(11):1413--1457, 2004.

\bibitem{daubechies2007projgrad}
Ingrid Daubechies, Massimo Fornasier, and Ignace Loris.
\newblock Accelerated projected gradient method for linear inverse problems
  with sparsity constraints.
\newblock To appear in \textit{Journal of Fourier Analysis and Applications},
  2008.

\bibitem{donoho2003sparse}
David~L. Donoho and Michael Elad.
\newblock Optimally-sparse representation in general (non-orthogonal)
  dictionaries via $\ell^1$ minimization.
\newblock {\em Proceedings of the National Academy of Sciences},
  100:2197--2202, 2003.

\bibitem{engl1996inverseproblems}
Heinz~W. Engl, Martin Hanke, and Andreas Neubauer.
\newblock {\em Regularization of Inverse Problems}, volume 375 of {\em
  Mathematics and its Applications}.
\newblock Kluwer Academic Publishers Group, Dordrecht, 2000.

\bibitem{fornasier2008jointsparsity}
Massimo Fornasier and Holger Rauhut.
\newblock Recovery algorithms for vector valued data with joint sparsity
  constraints.
\newblock {\em SIAM Journal on Numerical Analysis}, 46(2):577--613, 2008.

\bibitem{griesse2008ssnsparsity}
Roland Griesse and Dirk~A. Lorenz.
\newblock A semismooth {N}ewton method for {T}ikhonov functionals with sparsity
  constraints.
\newblock {\em Inverse Problems}, 24:035007 (19pp), 2008.

\bibitem{hofmann2007convtikban}
Bernd Hofmann, Barbara Kaltenbacher, Christiane Poeschl, and Otmar Scherzer.
\newblock A convergence rates result for {T}ikhonov regularization in {B}anach
  spaces with non-smooth operators.
\newblock {\em Inverse Problems}, 23(3):987--1010, 2007.

\bibitem{lorenz2007softhard}
Dirk~A. Lorenz.
\newblock Non-convex variational denoising of images: Interpolation between
  hard and soft wavelet shrinkage.
\newblock {\em Current Development in Theory and Application of Wavelets},
  1(1):31--56, 2007.

\bibitem{ramlau2005tikhreplace}
Ronny Ramlau and Gerd Teschke.
\newblock Tikhonov replacement functionals for iteratively solving nonlinear
  operator equations.
\newblock {\em Inverse Problems}, 21(5):1571--1592, 2005.

\bibitem{ramlau2006tikhproject}
Ronny Ramlau and Gerd Teschke.
\newblock A {T}ikhonov-based projection iteration for nonlinear ill-posed
  problems with sparsity constraints.
\newblock {\em Numerische Mathematik}, 104(2):177--203, 2006.

\bibitem{resmerita2005regbanspaces}
Elena Resmerita.
\newblock Regularization of ill-posed problems in {B}anach spaces: convergence
  rates.
\newblock {\em Inverse Problems}, 21(4):1303--1314, 2005.

\bibitem{resmerita2006nonquadreg}
Elena Resmerita and Otmar Scherzer.
\newblock Error estimates for non-quadratic regularization and the relation to
  enhancement.
\newblock {\em Inverse Problems}, 22(3):801--814, 2006.

\bibitem{schoepfer2006illposedbanach}
Frank Sch{\"o}pfer, Alfred~K. Louis, and Thomas Schuster.
\newblock Nonlinear iterative methods for linear ill-posed problems in {B}anach
  spaces.
\newblock {\em Inverse Problems}, 22:311--329, 2006.

\bibitem{schuster2008tikbanach}
Thomas Schuster, Peter Maass, Thomas Bonesky, Kamil~S. Kazimierski, and Frank
  Sch\"opfer.
\newblock Minimization of {T}ikhonov functionals in {B}anach spaces.
\newblock {\em Abstract and Applied Analysis}, 2008:Article ID 192679, 19
  pages, 2008.

\bibitem{xu1991charinequalities}
Zong~Ben Xu and Gary~F. Roach.
\newblock Characteristic inequalities of uniformly convex and uniformly smooth
  {Banach} spaces.
\newblock {\em Journal of Mathematical Analysis and Applications},
  157(1):189--210, 1991.

\end{thebibliography}

\end{document}